\documentstyle[fleqn,12pt]{article}
\begin{document}
\pagenumbering{arabic}
\setcounter{page}{1}
\pagestyle{plain}
\baselineskip=16pt

\thispagestyle{empty}
\begin{center}
{\Large\bf On the Differential Geometry of $GL_q(1\vert 1)$ } 
\end{center}

\vspace{1cm}
\noindent
Salih Celik$^*$ \footnote{{\bf E-mail}: scelik@fened.msu.edu.tr} 
and 
Sultan A. Celik$^+$ \footnote{{\bf E-mail}: celik@yildiz.edu.tr} 

\noindent
$^*$ {\footnotesize Mimar Sinan University, Department of Mathematics, 
80690 Besiktas, Istanbul, TURKEY.}\\
$^+$ {\footnotesize Yildiz Technical University, Department of Mathematics, 
80270 Sisli, Istanbul, TURKEY. }

\vspace{1.5cm}
{\bf Abstract.} 
The differential calculus on the quantum supergroup GL$_q(1\vert 1)$ was 
introduced by Schmidke {\it et al}. (1990 {\it Z. Phys. C} {\bf 48} 249). 
We construct a differential calculus on the quantum supergroup 
GL$_q(1\vert 1)$ in a different way and we obtain its quantum superalgebra. 
The main structures are derived without an R-matrix. It is seen that the 
found results can be written with help of a matrix $\hat{R}$

\vfill\eject
\noindent
{\bf 1 Introduction}

\noindent
During the past few years, the theory of quantum groups [1] has been an 
important branch of mathematical physics and a new branch of mathematics. 
A quantum group is not a group in ordinary sense of this word. But a quantum 
group is somehow related to a group in the sense that it is a deformation 
of a certain structure reflecting the group properties. In other words, a 
quantum group coincides with the group for particular values of the 
deformation parameter. Quantum (super) groups present the examples of 
(graded) Hopf algebras. They have found application in as diverse areas 
of physics and mathematics as conformal field theory, statistical mechanics, 
nonlinear integrable models, knot theory and solutions of Yang-Baxter 
equations [2] (and references therein). Many of the remarkable properties of 
Matrix Theory appear to be closely connected to the ideals of noncommutative 
geometry [3]. More recently it has been suggested that the zero branes in 
M-theory should be identified with supercoordinates in noncommutative 
geometry [4]. 

Quantum (super) groups can be realized on a quantum (super) space in which 
coordinates are noncommuting [5]. Recently the differential calculus on 
noncommutative space has been intensively studied both by mathematicians 
and mathematical physicists. There is much activity in differential geometry 
on quantum groups. A noncommutative differential calculus on the quantum 
groups has been developed by Woronowicz [6] along the lines of the general 
ideas of Connes [3]. Wess and Zumino [7] have reformulated the 
general theory in an abstract way. A few other methods to construct 
a noncommutative differential geometry on a quantum group have been 
proposed and discussed by several authors (e.g. Refs 8). 

In Ref. 9 a right-invariant differential calculus on the 
quantum supergroup $GL_q(1\vert 1)$ has been constructed 
and it has been shown that the quantum Lie algebra generators 
satisfy the {\it undeformed} Lie superalgebra. 
In this paper we present a differential calculus on the 
quantum supergroup $GL_q(1\vert 1)$ in a different way. This differential 
structure turns out to be a differential (graded) Hopf 
algebra. Although all of the commutation relations among 
the one-forms of Ref. 9 are classical, i.e. 'undeformed' 
in the present work they are $q$-deformed. It is also cited 
that the obtained relations can be written with help of a 
matrix $\hat R$. 

In this work, we shall use greek letters to denote fermionic (odd) and 
latin letters to denote bosonic (even) variables. 

\noindent
{\bf 2 Formalities}

\noindent
Elementary properties of the quantum supergroup $GL_q(1\vert 1)$ are 
described in [10,11]. We state briefly the properties we are going to 
need in this work. 

The quantum supergroup $GL_q(1\vert 1)$ is defined by, as a group element, 
the matrices of the form 
$$T = \left(\matrix{ a & \beta \cr \gamma & d \cr} \right) \eqno(1)$$ 
where the matrix elements satisfy the following commutation 
relations [10,11] 
$$ a \beta = q \beta a \qquad d \beta = q \beta d $$
$$ a \gamma = q \gamma a \qquad d \gamma = q \gamma d \eqno(2)$$
$$ \beta \gamma + \gamma \beta = 0 \qquad \beta^2 = 0 = \gamma^2 $$
$$ a d = d a + (q - q^{-1}) \gamma \beta. $$
Let us denote the algebra generated by the elements $a$, 
$\beta$, $\gamma$, $d$ with the relations (2) by ${\cal A}$. 
We know that the algebra 
${\cal A}$ is a (graded) Hopf algebra with the following structure: 

{\bf (1)} The usual coproduct 
$$\Delta : {\cal A} \longrightarrow {\cal A} \otimes {\cal A} \qquad 
   \Delta(T) = T \dot{\otimes} T  \eqno(3)$$

{\bf (2)} the counit 
$$ \varepsilon : {\cal A} \longrightarrow {\cal C} \qquad 
    \varepsilon(T) = I \eqno(4)$$

{\bf (3)} the coinverse (antipode) 
$$S : {\cal A} \longrightarrow {\cal A} \qquad 
  S(T) = 
  \left(\matrix{ 
    a^{-1} + a^{-1} \beta d^{-1} \gamma a^{-1} & - a^{-1} \beta d^{-1} \cr 
    - d^{-1} \gamma a^{-1} & d^{-1} + d^{-1} \gamma a^{-1} \beta d^{-1} \cr }
     \right).  \eqno(5)$$
It is not difficult to verify the following properties of the 
co-structures: 
$$(\Delta \otimes \mbox{id}) \circ \Delta = 
  (\mbox{id} \otimes \Delta) \circ \Delta \eqno(6)$$
$$\mu \circ (\varepsilon \otimes \mbox{id}) \circ \Delta 
  = \mu' \circ (\mbox{id} \otimes \varepsilon) \circ \Delta \eqno(7)$$
$$m \circ (S \otimes \mbox{id}) \circ \Delta = \varepsilon 
  = m \circ (\mbox{id} \otimes S) \circ \Delta \eqno(8)$$
where id denotes the identity mapping, 
$$\mu : {\cal C} \otimes {\cal A} \longrightarrow {\cal A} \qquad 
  \mu' : {\cal A} \otimes {\cal C} \longrightarrow {\cal A} $$
are the canonical isomorphisms, defined by 
$$\mu(k \otimes a) = ka = \mu'(a \otimes k) \quad \forall a \in {\cal A} 
  \quad \forall k \in {\cal C} \eqno(9)$$
and $m$ is the multiplication map 
$$m : {\cal A} \otimes {\cal A} \longrightarrow {\cal A} \qquad 
  m(a \otimes b) = ab. $$
The multiplication in ${\cal A} \otimes {\cal A}$ follows the rule 
$$(A \otimes B) (C \otimes D) = (-1)^{p(B) p(C)} AC \otimes BD \eqno(10)$$
where $p(X)$ is the $z_2$-grade of $X$, i.e. $p(X) = 0$ for even variables 
and $p(X) = 1$ for odd variables. 

In the following section we shall build up the differential calculus on the 
quantum supergroup $GL_q(1\vert 1)$. We first note that the properties of 
the exterior differential. The exterior differential {\sf d} is an operator 
which gives the mapping from the generators of ${\cal A}$ to the 
differentials: 
$${\sf d} : u \longrightarrow {\sf d} u \qquad u \in \{a,\beta,\gamma,d\}. 
  \eqno(11)$$
We demand that the exterior differential has to satisfy two properties: the 
nilpotency 
$${\sf d}^2 = 0 \eqno(12)$$
and the graded Leibniz rule 
$${\sf d}(fg) = ({\sf d} f) g + (-1)^{p(f)} f ({\sf d} g). \eqno(13)$$

We now introduce the algebra $\hat{\cal A}$ generated by the 
matrix elements $\alpha$, $b$, $c$, $\delta$ of $\hat {T}$, 
$$ \hat{T} = \left(\matrix{ \alpha & b \cr c & \delta \cr} \right) \eqno(14)$$
where the matrix elements $\alpha$, $b$, $c$, $\delta$ satisfy the 
commutation relations [12] 
$$ \alpha b = q^{-1} b \alpha \qquad \alpha c = q^{-1} c \alpha $$
$$ \delta b = q^{-1} b \delta \qquad \delta c = q^{-1} c \delta $$
$$\alpha \delta + \delta \alpha = 0 \qquad \alpha^2 = 0 = \delta^2 
   \eqno(15)$$
$$ bc = cb + (q - q^{-1}) \delta \alpha. $$
We shall use these relations in the following section. 

\vfill\eject\noindent
{\bf 3 Differential Geometric Structure of $GL_q(1\vert 1)$}

\noindent
We have seen, in the previous section, that ${\cal A}$ is an associative 
algebra (essentially a graded Hopf algebra) generated by the matrix 
elements of (1) with the relations (2). A differential algebra on 
${\cal A}$ is a $z_2$-graded associative algebra $\Omega$ equipped 
with a linear operator {\sf d} given (11)-(13). Furthermore the algebra 
$\Omega$ has to be generated by $\Omega^0 \cup {\sf d} \Omega^0$, where 
$\Omega^0$ is isomorphic to ${\cal A}$. 

Firstly, we observe that the matrix elements of $\hat{T}$ are given by (14) 
just as the differentials of the matrix elements of $T$ are given by (1). 
Thus, as in the considerations for the quantum planes in [7], we can 
identify, at least formally, 
$$\hat{T} = {\sf d} T. \eqno(16)$$
So, let us rewrite the relations (15) in the form 
$${\sf d} a {\sf d} \beta = q^{-1} {\sf d} \beta {\sf d} a \qquad 
   {\sf d} d {\sf d} \beta = q^{-1} {\sf d} \beta {\sf d} d $$
$${\sf d} a {\sf d} \gamma = q^{-1} {\sf d} \gamma {\sf d} a \qquad 
  {\sf d} d {\sf d} \gamma = q^{-1} {\sf d} \gamma  {\sf d} d \eqno(17)$$
$${\sf d}a {\sf d} d = - {\sf d} d {\sf d} a \qquad 
   ({\sf d} a)^2 = 0 =  ({\sf d} d)^2 $$
$$ {\sf d}\beta {\sf d} \gamma = {\sf d} \gamma {\sf d} \beta + 
   (q - q^{-1}) {\sf d} d {\sf d} a. $$

We now denote by ${\cal A}_{a\beta}$ the algebra generated by the elements 
$a$ and $\beta$ with the relations 
$$a \beta = q \beta a \qquad \beta^2 = 0. \eqno(18)$$
A possible set of commutation relations between the generators of 
${\cal A}_{a\beta} \simeq \Omega^0_{a\beta}$ and 
$\hat{\cal A}_{{\sf d}a {\sf d}\beta} \simeq 
{\sf d} \Omega^0_{{\sf d}a {\sf d}\beta}$ has 
of the form 
$$a~ {\sf d}a = A {\sf d}a~ a $$
$$a~ {\sf d}\beta = F_{11} {\sf d}\beta~ a + F_{12} {\sf d}a~ \beta 
  \eqno(19)$$
$$\beta~ {\sf d}a = F_{21} {\sf d}a~ \beta + F_{22} {\sf d}\beta~ a $$
$$\beta~ {\sf d}\beta = B {\sf d}\beta~ \beta. $$
The coefficients $A, B$ and $F_{ij}$ will be determined in term of complex 
deformation parameter $q$. To find them we shall use the consistency of 
calculus. From the consistency conditions 
$${\sf d}(a \beta - q \beta a) = 0 \qquad {\sf d}(\beta^2) = 0 \eqno(20)$$
and 
$$(a \beta - q \beta a) {\sf d}a = 0 \qquad 
  (a \beta - q \beta a) {\sf d}\beta = 0 \eqno(21)$$
we find 
$$F_{11} + q F_{22} = q \qquad F_{12} + q F_{21} = - 1 \qquad B = 1 \eqno(22)$$
and 
$$F_{12} F_{22} = 0 = (F_{11} - q A)F_{22}. \eqno(23)$$
In fact (22) and (23) are, in disguise, the linear and quadratic consistency 
conditions similar to the ones discussed in full generality for quantum 
planes in [7]. Equation (23) admits two solutions. If we choose $F_{22}  = 0$ 
we are led to the following commutation relations 
$$ a ~{\sf d}a = q^2 {\sf d}a~ a $$
$$ a~ {\sf d}\beta = q {\sf d}\beta~ a + (q^2 - 1) {\sf d}a~ \beta\eqno(24a)$$
$$\beta~ {\sf d}a = - q {\sf d}a~ \beta $$
$$ \beta~ {\sf d}\beta = {\sf d}\beta~ \beta $$
where $A$ equal to $q^2$ since this leads to the standard R-matrix (eq. (51)). 

We denote by 
$\Omega_{a\beta} = {\cal A}_{a\beta} \cup 
  {\sf d}{\cal A}_{{\sf d}a {\sf d}\beta}$, 
the algebra generated by the generators $a$, $\beta$ (of ${\cal A}$) and 
the generators ${\sf d}a$, ${\sf d}\beta$ (of ${\sf d}{\cal A})$ with the 
relations (18), a parity of (17) and (24a). 

Continuing in this way, we can construct the differential algebra 
$\Omega_{a\gamma}$, $\Omega_{\beta\gamma}$, etc. The final result is 
given by 
$$ a ~{\sf d}\gamma = q {\sf d}\gamma~ a + (q^2 - 1) {\sf d}a~ \gamma$$
$$a~ {\sf d}d = {\sf d}d~ a + (q - q^{-1}) [{\sf d}\gamma~ \beta 
  - {\sf d}\beta~ \gamma + (q - q^{-1}) {\sf d}a~ d]$$
$$ \beta~ {\sf d}\gamma = {\sf d}\gamma~ \beta + (q - q^{-1}){\sf d}a~ d $$
$$ \beta ~{\sf d}d = - q^{-1} {\sf d}d~ \beta + (1 - q^{-2}) {\sf d}\beta~ d$$
$$\gamma~ {\sf d}a = - q {\sf d}a~ \gamma \qquad 
  \gamma~ {\sf d}\beta = {\sf d}\beta~ \gamma - (q - q^{-1}){\sf d}a~ d 
  \eqno(24b)$$
$$ \gamma ~{\sf d}\gamma = {\sf d}\gamma~ \gamma \qquad 
 \gamma~ {\sf d}d = - q^{-1} {\sf d}d~ \gamma + (1 - q^{-2}) {\sf d}\gamma~ d$$
$$ d~ {\sf d}a = {\sf d}a~ d \qquad d~ {\sf d}\beta = q^{-1} {\sf d}\beta~ d$$
$$ d~ {\sf d}\gamma = q^{-1} {\sf d}\gamma~ d \qquad 
   d~ {\sf d}d = q^{-2} {\sf d}d~ d.$$ 
Thus we have constructed the differential algebra of the algebra 
${\cal A}$ generated by the matrix elements of any matrix in $GL_q(1\vert 1)$. 
It is not difficult to check that the action of {\sf d} on (24) is consistent. 

We note that the differential algebra $\Omega = {\cal A} \cup {\sf d}{\cal A}$ 
is a (graded) Hopf algebra with the following co-structure [13]: the 
coproduct is given by 
$$ \hat{\Delta} : \Omega \longrightarrow \Omega \otimes \Omega 
   \qquad 
   \hat{\Delta}(\hat{T}) = \hat{T} \dot{\otimes} T + 
   (-1)^{p(T)} T \dot{\otimes} \hat{T}. 
   \eqno(25)$$
Explicity, 
$$\hat{\Delta}({\sf d}a) = {\sf d}a \otimes a + {\sf d}\beta \otimes \gamma 
  + a \otimes {\sf d} a - \beta \otimes {\sf d}\gamma $$ 
$$\hat{\Delta}({\sf d}\beta) = {\sf d}\beta \otimes d + {\sf d}a \otimes \beta 
  + a \otimes {\sf d}\beta - \beta \otimes {\sf d}d \eqno(26)$$ 
$$\hat{\Delta}({\sf d}\gamma) = {\sf d}\gamma \otimes a + 
  {\sf d}d \otimes \gamma - \gamma \otimes {\sf d}a + d \otimes {\sf d}\gamma$$
$$\hat{\Delta}({\sf d}d) = {\sf d}d \otimes d + {\sf d}\gamma \otimes \beta - 
   \gamma \otimes {\sf d}\beta + d \otimes {\sf d}d. $$ 
The counit is given by 
$$\hat{\varepsilon} : \Omega \longrightarrow {\cal C}  \qquad 
  \hat{\varepsilon}(\hat{T}) = 0 \eqno(27)$$
and the coinverse 
$$\hat{S} : \Omega  \longrightarrow \Omega  \qquad 
  \hat{S}(\hat{T}) = - (-1)^{p(T^{-1})} T^{-1} \hat{T} T^{-1}. \eqno(28)$$
The central element is 
$$\hat{D} = b c^{-1} - \alpha c^{-1} \delta c^{-1} \eqno(29)$$
where $\hat{T} = {\sf d}T$, i.e. ${\hat D}$ commutes with the generators of 
${\cal A}$ and also with the generators of $\hat{\cal A}$. 

It is not difficult to verify that the maps $\hat{\Delta}$ and 
$\hat{\varepsilon}$ are both algebra homomorphisms and $\hat{S}$ is an 
algebra anti-homomorphism. The three maps satisfy also the properties 
(6)-(8), and they preserve the relations (24) provided that the 
actions on the generators of ${\cal A}$ of $\hat{\Delta}$, 
$\hat{\varepsilon}$ and $\hat{S}$ are the same as (3)-(5). 
                    
\noindent
{\bf 4 The Cartan-Maurer Forms in $\Omega$}

\noindent
In analogy with the right invariant one-forms defined on a Lie group in 
classical differential geometry, one can construct the matrix valued 
one-form $W$ where 
$$W = {\sf d}T T^{-1}. \eqno(30)$$
If we set 
$$T^{-1} = \left(\matrix{ 
    A & B \cr 
    C & D \cr} \right) \eqno(31)$$
as the superinverse of $T \in GL_q(1\vert 1)$, we write the matrix elements 
(one-forms) of $W$ 
$$w_1 = {\sf d}a A + {\sf d}\beta C \qquad u = {\sf d}a B + {\sf d}\beta D $$
$$w_2 = {\sf d}d D + {\sf d}\gamma B \qquad 
  v = {\sf d}\gamma A + {\sf d}d C. \eqno(32)$$
Firstly, we obtain the commutation relations between the matrix elements of 
$T$ and $T^{-1}$ as follows: 
$$a A = q^2 A a + 1 - q^2 \qquad d A = A d $$
$$a D = D a \qquad d D = q^2 D d + 1 - q^2 $$
$$a B = q B a \qquad d B = q B d $$
$$a C = q C a \qquad d C = q C d $$
$$\beta A = q A \beta \qquad \gamma A = q A \gamma \eqno(33)$$
$$\beta D = q D \beta \qquad \gamma D = q D \gamma$$
$$\beta B = B \beta \qquad \gamma B = - q^2 B \gamma$$
$$\beta C = - q^2 C \beta \qquad \gamma C = C \gamma.$$

Using these relations, we now find the commutation relations of the matrix 
entries of $T$ with those of $W$ : 
$$a w_1 = q^2 w_1 a \qquad d w_2 = q^{-2} w_2 d + (q^{-2} - 1) v \beta$$
$$a w_2 = w_2 a + (q^{-2} - 1) u \gamma + (q - q^{-1})^2 w_1 a \qquad 
   d w_1 = w_1 d$$
$$a u = q u a \qquad du = q^{-1} ud - (q - q^{-1}) w_1 \beta $$
$$a v = q v a + (q - q^{-1}) w_1 \gamma \qquad 
  d v = q^{-1} v d \eqno(34)$$
$$\beta w_1 = - q^2 w_1 \beta \qquad \gamma w_2 = - q^{-2} w_2 \gamma + 
   (1 - q^{-2}) v a$$
$$\beta w_2 = - w_2 \beta + (1 - q^{-2}) u d - (q - q^{-1})^2 w_1 \beta \qquad 
   \gamma w_1 = - w_1 \gamma$$
$$\beta u = q u \beta \qquad \gamma u = q^{-1} u \gamma - 
   (q - q^{-1}) w_1 a$$
$$\beta v = q v \beta + (q - q^{-1}) w_1 d \qquad 
   \gamma v = q^{-1} v \gamma.$$
These relations are different from those of [9] excepting a few of them. 

To obtain commutation relations among the Cartan-Maurer one-forms we shall 
use the commutation relations between the matrix elements of $T^{-1}$ and 
the differentials of the matrix elements of $T$ which are given in the 
following 
$$A {\sf d}a = q^{-2} {\sf d}a A \qquad A {\sf d}d = {\sf d}d A $$
$$A {\sf d}\beta = q^{-1} {\sf d}\beta A \qquad 
  A {\sf d}\gamma = q^{-1} {\sf d}\gamma A $$
$$D {\sf d}a = {\sf d}a D \qquad 
  D {\sf d}d = q^2 {\sf d}d D + (q^2 - 1) ({\sf d}\gamma B - {\sf d}\beta C) 
               - (q - q^{-1})^2 {\sf d}a A$$
$$D {\sf d}\beta = q {\sf d}\beta D + (q - q^{-1}) {\sf d}a B \qquad 
  D {\sf d}\gamma = q {\sf d}\gamma D + (q - q^{-1}) {\sf d}a C $$
$$B {\sf d}a = - q^{-1} {\sf d}a B \qquad 
  B {\sf d}\gamma = {\sf d}\gamma B + (q^{-2} - 1) {\sf d}a A \eqno(35)$$
$$B {\sf d}\beta = {\sf d}\beta B \qquad 
  B {\sf d}d = - q {\sf d}d B - (q - q^{-1}) {\sf d}\beta A$$
$$C {\sf d}a = - q^{-1} {\sf d}a C \qquad 
  C {\sf d}\beta = {\sf d}\beta C + (1 - q^{-2}) {\sf d}a A $$
$$C {\sf d}\gamma = {\sf d}\gamma C \qquad 
  C {\sf d}d = - q {\sf d}d C - (q - q^{-1}) {\sf d}\gamma A.$$

We now obtain the commutation relations of the Cartan-Maurer forms, using 
(35) and (17), as follows: 
$$w_1 u = u w_1 \qquad u w_2 = q^2 w_2 u + (1 - q^2) w_1 u $$
$$w_1 v = v w_1 \qquad w_2 v = q^2 v w_2 + (1 - q^2) v w_1 $$
$$w_1^2 = 0 \qquad w_2^2 = (1 - q^2) v u \eqno(36)$$
$$u v = q^2 v u \qquad w_1 w_2 + w_2 w_1 = (1 - q^2) v u. $$
Again, one can check that the action of {\sf d} on (34) and also on (36) is 
consistent. 

Note that although the Cartan-Maurer one-forms of Ref. 9 satisfy 
$q$-independent commutation relations, i.e. they are identical to the 
classical ones, in the present work they satisfy the $q$-dependent relations. 

\noindent
{\bf 5 Quantum Superalgebra}

\noindent
In this section we shall obtain the deformed superalgebra of the Lie 
algebra generators. 
Firstly, we write, from (30), the Cartan-Maurer forms as 
$${\sf d}a = w_1 a + u \gamma \qquad {\sf d}\beta = w_1 \beta + u d$$
$${\sf d}d = w_2 d + v \beta \qquad {\sf d}\gamma = w_2 \gamma + v a. 
\eqno(37)$$
We can write ${\sf d}W$ in the form 
$${\sf d}W = \sigma_3 W \sigma_3 W \qquad 
  \sigma_3 = \left(\matrix{1 & 0 \cr 0 & - 1\cr}\right) \eqno(38)$$
In terms of the two-forms, these become 
$${\sf d}w_1 = w_1^2 - u v \qquad {\sf d}u = w_1 u - u w_2$$
$${\sf d}w_2 = w_2^2 - v u \qquad {\sf d}v = w_2 v - v w_1. \eqno(39)$$
We can now write down the Cartan-Maurer equations in our case 
$$\left(\matrix { {\sf d}w_1 & {\sf d}u   \cr 
                  {\sf d}v   & {\sf d}w_2 \cr}\right) 
  = \left(\matrix{ - u v          & q^2 (w_1 - w_2) u \cr 
                  - (w_1 - w_2) v & - u v \cr}\right). \eqno(40)$$

The commutation relations of the Cartan-Maurer forms allow us to 
construct the algebra of the generators. To obtain the quantum 
superalgebra we write down the exterior differential in the form 
$${\sf d} = w_1 T_1 + w_2 T_2 + u \nabla_+ + v \nabla_-. \eqno(41)$$
Considering an arbitrary function $f$ of the group parameters and using the 
nilpotency of the exterior differential {\sf d} one has 
$$({\sf d}w_i) T_if + ({\sf d}u_i) \nabla_if = 
  w_i {\sf d}T_if - u_i {\sf d}\nabla_i \eqno(42)$$
where 
$$w_i \in \{w_1,w_2\} \qquad u_i \in \{u,v\} \qquad 
  \nabla_i \in \{\nabla_+, \nabla_-\}.$$
With the Cartan-Maurer equations we find the following commutation 
relations for the quantum superalgebra: 
$$[T_1,\nabla_+] = - q^2 \nabla_+ + (q^2 - 1) T_2 \nabla_+ $$
$$[T_2,\nabla_+] = q^2 \nabla_+ + (1 - q^2) T_2 \nabla_+ $$
$$[T_1,\nabla_-] = q^2 \nabla_- + (1 - q^2) \nabla_- T_2$$
$$[T_2,\nabla_-] = - q^2 \nabla_- + (q^2 - 1) \nabla_- T_2 \eqno(43)$$
$$[T_1,T_2] = 0 \qquad \nabla_{\pm}^2 = 0$$
$$\nabla_- \nabla_+ + q^{-2} \nabla_+ \nabla_- = T_1 + T_2 + 
  (q^{-2} - 1) (T_2^2 + T_1T_2) $$
or with new generators $X = T_1 + T_2$ and $Y = T_1 - T_2$, 
$$[X,\nabla_{\pm}] = 0 \qquad [X,Y] = 0 \qquad \nabla_{\pm}^2 = 0$$
$$[Y,\nabla_+] = - 2 q^2 \nabla_+ + (q^2 - 1) (X - Y)\nabla_+$$
$$[Y,\nabla_-] = 2 q^2 \nabla_- + (1 - q^2) \nabla_- (X - Y) \eqno(44)$$
$$\nabla_+ \nabla_- + q^2 \nabla_- \nabla_+ = q^2 X + (1 - q^2) (X^2 - XY). $$

We also note that the commutation relations (43) of the superalgebra 
generators should be consistent with monomials of the group parameters. To 
proceed, we must calculate the actions of the Leibniz  rule by comparing 
the elements which lie together with each other from the Cartan-Maurer froms: 
$$T_1 a  = a + q^2 aT_1 + (q - q^{-1})^2 a T_2 + (q - q^{-1}) \gamma \nabla_-$$
$$T_1 \beta  = \beta + q^2 \beta T_1 - (q - q^{-1})^2 \beta T_2 - 
  (q - q^{-1}) d \nabla_-$$
$$T_1 \gamma  = \gamma T_1 - (q - q^{-1}) a \nabla_+ \qquad 
  T_1 d  = d T_1 - (q - q^{-1}) \beta \nabla_+ $$
$$T_2 a  = a T_2 \qquad T_2 d = d + q^{-2} d T_2 $$
$$T_2 \beta  = \beta T_2 \qquad 
  T_2 \gamma  = \gamma + q^{-2} \gamma T_2 \eqno(45)$$
$$\nabla_+ a = \gamma + q a \nabla_+ + (q^{-2} - 1) \gamma T_2 \qquad 
  \nabla_+ d = q^{-1} d \nabla_+ $$
$$\nabla_+ \beta = d - q \beta \nabla_+ + (q^{-2} - 1) d T_2 \qquad 
  \nabla_+ \gamma = - q^{-1} \gamma \nabla_+ $$
$$\nabla_- a = q a \nabla_- \qquad 
  \nabla_- d = \beta + q^{-1} d \nabla_- + (q^{-2} - 1) \beta T_2 $$
$$\nabla_- \beta = - q \beta \nabla_- \qquad 
  \nabla_- \gamma = a - q^{-1} \gamma \nabla_- + (q^{-2} - 1) a T_2. $$

Notice that these commutation relations must be consistent. In fact, it is 
easy to see that the nilpotency of $\nabla_\pm^2 = 0$ is consistent with 
$$\nabla_\pm^2 a = q^2 a \nabla_\pm^2. $$
Also, the commutation relations of $T_1$, $T_2$ with $\nabla_-$ are consistent 
with 
$$\left([T_1,\nabla_-] - q^2 \nabla_- + (q^2 - 1)\nabla_- T_2\right)a = 
 q a \left([T_1,\nabla_-] - q^2 \nabla_- + (q^2 - 1)\nabla_- T_2\right)$$
$$\left(T_2 \nabla_- - q^2 \nabla_- T_2 + q^2 \nabla_- \right) a = 
 q a \left(T_2 \nabla_- - q^2 \nabla_- T_2 + q^2 \nabla_- \right).$$
Furthermore, 
$$(T_1 T_2 - T_2 T_1) a = (q^{-3} - q^{-1}) \gamma 
  \left(T_2 \nabla_- - q^2 \nabla_- T_2 + q^2 \nabla_- \right)$$
and 

$\left(\nabla_- \nabla_+ + q^{-2} \nabla_+ \nabla_- - T_1 - T_2 - 
  (q^{-2} - 1) (T^2_2 + T_1 T_2)\right) a $

$\hspace*{1.9cm}= q^2 a 
\left(\nabla_- \nabla_+ + q^{-2} \nabla_+ \nabla_- - T_1 - T_2 - 
  (q^{-2} - 1) (T^2_2 + T_1 T_2)\right).$ 

\noindent
Similarly, one can find the other relations. 

\noindent
{\bf 6 R-matrix approach}

\noindent
In this section we wish to obtain the relations (17), (24), (34) and (36) 
with the help of a matrix $\hat{R}$ that acts on the square tensor space 
of the supergroup. Of course, the matrix $\hat R$ is a solution of the 
quantum (graded) braided group equation. 

We first consider the quantum superplane and its dual [5]. The quantum 
superplane $A_q$ is generated by coordinates $x$ and $\theta$, and the 
commutation rules 
$$x \theta = q \theta x \qquad \theta^2 = 0. \eqno(46)$$
The quantum (dual) superplane $A^*_q$ is generated by coordinates $\varphi$ 
and $y$, and the commutation rules 
$$\varphi^2 = 0 \qquad \varphi y = q^{-1} y \varphi. \eqno(47)$$
We demand that relations (46), (47) are preserved under the action of $T$, 
as a linear transformation, on the quantum superplane and its dual: 
$$T : A_q \longrightarrow A_q \qquad T : A^*_q \longrightarrow A^*_q. 
  \eqno(48)$$
Let $X = (x, \theta)^t$ and $\hat{X} = (\varphi, y)^t$. Then, as a 
consequence of (48) the points $TX$ and $T\hat{X}$ should belong to 
$A_q$ and $A_q^*$, respectively, which give the relations (2). 

Similarly, let us consider linear transformations $\hat T$ with the 
following properties 
$$\hat{T} : A_q \longrightarrow A^*_q \qquad 
  \hat{T} : A^*_q \longrightarrow A_q.   \eqno(49)$$
Then the points $\hat{T}X$ and $\hat{T}\hat{X}$ should belong to $A_q^*$ 
and $A_q$, respectively. This case is equivalent to (17). 

Note that the relations (46) can be written as follows 
$$X \otimes X = q^{-1} \hat{R} X \otimes X \eqno(50)$$
where 
$$ \hat{R} = \left(\matrix{ 
   q &    0       & 0 & 0 \cr
   0 & q - q^{-1} & 1 & 0 \cr
   0 &    1       & 0 & 0 \cr
   0 &    0       & 0 & -q^{-1} \cr} \right). \eqno(51)$$
We can also write mixed relations between the component of $X$ and $\hat X$ 
as follows: 
$$(-1)^{p(X)} X \otimes \hat{X} = q \hat{R} \hat{X} \otimes X \eqno(52)$$
where $\hat{X} = {\sf d}X$. 

Using (48) together with (50) and (52), we now derive anew the quantum 
supergroup relations (2) from the equation [9] 
$$\hat {R} T_1 T_2 = T_1 T_2 \hat {R} \eqno(53)$$
where, in usual grading tensor notation, $T_1 = T \otimes I$ and 
$T_2 = I \otimes T$. Similarly using (52), we obtain the following equation 
$$(-1)^{p(T_1)} T_1 \hat {T}_2 = \hat {R} \hat {T}_1 T_2 \hat {R} \eqno(54)$$
which is equivalent to the relations (24). The equation 
$$\hat {T}'_1 \hat {T}_2 = 
  (-1)^{p(\hat {T}_1)} \hat {R} \hat {T}_1 \hat {T}_2\hat {R} 
  \qquad \hat {T}'_1 = {\sf d}((-1)^{p(T_1)} T_1) \eqno(55)$$
gives the relations (17). Taking $\hat {T} = W T$ and using (55) one has 
$$(-1)^{p(T_1)} T_1 W_2 = \hat {R} W_1 \hat {R} T_1 \eqno(56)$$
which gives the relations (34). Finally, from (55) we find that 
$$(-1)^{p(W_1)} W_1 \hat {R} W_1 \hat {R}^{-1} + 
  (-1)^{p(W_1)} \hat {R} W_1 \hat {R} W_1 = 0. \eqno(57)$$
This equation is equivalent to (36). 

\noindent
{\bf 7 Discussion} 

\noindent
The present paper may be considered as an alternative to the approach proposed 
earlier by Schmidke {\it et al} [9]. They have consructed a right-invariant 
differential calculus on the quantum supergroup $GL_q(1\vert 1)$ and showed 
that the quantum superalgebra generators satisfy the undeformed Lie 
superalgebra. Their starting point is to evaluate the $q$-commutation 
relations of the group parameters [the matrix elements of the matrix in 
$GL_q(1\vert 1)$] with the Cartan-Maurer forms, directly. They have assumed 
that "the right action of the group suggests that $a$ and $\beta$, and, 
$\gamma$ and $d$ [the matrix elements in (1)] satisfy the same relations 
as the forms". The starting point of the present paper, however, is to 
evaluate the $q$-commutation relations of the matrix elements with their 
differentials. Later, using these relations the $q$-commutation relations of 
the matrix elements with the Cartan-Maurer forms are obtained without any 
further assumptions. In the work of [9], the commutation relations of the 
Cartan-Maurer forms are obtained by using the $q$-commutation relations 
of the matrix elements with the Cartan-Maurer forms, i.e., to obtain the 
desired commutation relations they have applied the exterior differential 
$\delta$ on the relations of the matrix elements with the Cartan-Maurer forms. 
In the present work we also use the same approach. But although the 
Cartan-Maurer forms in [9] satisfy $q$-independent relations, in the present 
paper they satisfy $q$-dependent relations. Since in [9], the Cartan-Maurer 
forms satisfy undeformed [$q$-independent] relations, the quantum Lie algebra 
obtained is also undeformed. 

The work of [9] does not allow an R-matrix approach. In our work we have 
derived the $q$-commutation relations between the matrix elements and their 
differentials without considering an R-matrix at first. However we later 
show that these relations can also be derived using an R-matrix. This 
R-matrix turns out to be the same as used by Wess and Zumino [7] for the 
commutation relations between coordinates and their differentials for 
the $GL_q(2)$ invariant calculus on the $q$-plane.

\noindent
{\bf Acknowledgments}

\noindent
This work was supported in part by {\bf T. B. T. A. K.} the 
Turkish Scientific and Technical Research Council. 
We would like to express our deep gratitude to the referees for critical 
comments and suggestions on the maniscript.

\end{document}